\documentclass{amsart}

\usepackage[latin1]{inputenc}

\usepackage{amssymb}
\usepackage{amsmath}
\usepackage{amsthm}

\theoremstyle{plain}
\newtheorem{theorem}{Theorem}
\newtheorem{proposition}[theorem]{Proposition}

\theoremstyle{remark}
\newtheorem{remark}{Remark}

\newcommand{\vect}[1]{\ensuremath{\mathbf{#1}}} 
\newcommand{\minor}[3]{\ensuremath{{#1}_{{#2} \gets {#3}}}} 
\DeclareMathOperator{\ess}{ess} 
\DeclareMathOperator{\essl}{ess^{<}} 
\DeclareMathOperator{\gap}{gap} 

\DeclareMathOperator{\median}{median}

\begin{document}
\title[The arity gap of lattice polynomial functions]{The arity gap of polynomial functions over bounded distributive lattices}

\author{Miguel Couceiro}
\address[M. Couceiro]{Mathematics Research Unit \\
University of Luxembourg \\
6, rue Richard Coudenhove-Kalergi \\
L-1359 Luxembourg \\
Luxembourg}
\email{miguel.couceiro@uni.lu}
\author{Erkko Lehtonen}
\address[E. Lehtonen]{Computer Science and Communications Research Unit \\
University of Luxembourg \\
6, rue Richard Coudenhove-Kalergi \\
L-1359 Luxembourg \\
Luxembourg}
\email{erkko.lehtonen@uni.lu}

\begin{abstract}
Let $A$ and $B$ be arbitrary sets with at least two elements. 
The arity gap of a function $f\colon A^n\to B$ is the minimum decrease in its essential arity when essential arguments of $f$ are identified.
In this paper we study the arity gap of polynomial functions over bounded distributive lattices and present a complete classification of such functions in terms of their arity gap. To this extent, we present a characterization of the essential arguments of polynomial functions, which we then use to show that almost all lattice polynomial functions have arity gap 1, with the exception of truncated median functions, whose arity gap is 2.
\end{abstract}

\maketitle

\section{Introduction}
\label{sec:intro}
Current research in many-valued logic and computer science led to investigations in the theory of essential variables in several directions concerning, in particular, the distribution of values of functions whose variables are all essential (see, e.g., \cite{Davies, Salomaa, Yablonski}), the process of substituting variables for constants (see, e.g., \cite{Cimev79,Cimev86, Lupanov, Salomaa, Solovjev}) and  
the process of substituting variables for variables (see, e.g., \cite{CL, DK, Salomaa, Willard}).

The latter line of study goes back to the 1963 paper by Salomaa \cite{Salomaa} who considered the following problem: How is the number of essential variables of a given finite function $f$ affected when variables of $f$ are identified? The minimum decrease in the number of essential variables of $f$ when essential variables are identified is referred to as the \emph{arity gap} of $f$. Using a result by Salomaa \cite{Salomaa} concerning the substitution of variables for constants, it was shown that the arity gap of any function $f \colon A^n \to B$ is at most $\lvert A \rvert$ (see~\cite{CL}). Salomaa~\cite{Salomaa} provided examples of functions meeting this upper bound. Willard~\cite{Willard} improved this upper bound showing that the arity gap of $f$ is at most $2$ whenever $f$ has more than $\lvert A \rvert$ essential variables.

The upper bound $\lvert A \rvert$ on the arity gap of $f \colon A^n \to B$ tells us that, even though we cannot hope to fully classify functions according to their arity gap, such a complete classification may be achieved by imposing certain conditions on the functions considered.
One approach is to focus on finite functions, i.e., to require that $A$ is finite. In this setting, several efforts have been made in this direction (see, e.g., \cite{BK, GK, Willard}) which culminated in a complete classification of finite functions according to their arity gap (see \cite{CL2}).
Another approach is to focus on specific classes of functions while allowing arbitrary domains. 

In this paper, we take the latter approach and study the arity gap of polynomial functions over bounded distributive lattices.
In Section~\ref{sec:basic}, we recall the basic notions needed throughout the manuscript and present the classifications of functions with Boolean variables according to their arity gap as given in \cite{CL} and \cite{CL2}. 
In Section~\ref{sec:gap}, we focus on polynomial functions over arbritary bounded distributive lattices. We start by recalling canonical representations of polynomial functions over bounded distributive lattices as obtained by Goodstein \cite{Goo67}, and which we then use to describe the essential arguments of polynomial functions. The complete classification of polynomial functions according to their arity gap is given in Subsection~\ref{ssec:gapl}, 
Theorem \ref{GapPol}, which asserts that the only polynomial functions having arity gap 2 are the truncated median functions, i.e., functions of the form
\[
(a \vee {\median}) \wedge b.
\]
All other polynomial functions have arity gap 1.

\section{Basic notions and preliminary results}
\label{sec:basic}

Let $A$  and $B$ be arbitrary sets with at least two elements. 
By a \emph{$B$-valued function on $A$} we mean a mapping $f \colon A^n \to B$ for some positive integer $n$, called the \emph{arity} of $f$. If $B = A$, we refer to these functions as \emph{operations on $A$.} A typical example is the $i$-th $n$-ary \emph{projection}, that is, the mapping $(a_1, \ldots, a_n) \mapsto a_i$, denoted by $x_i^{(n)}$, or simply by $x_i$ when the arity is clear from the context.  
Operations on the two-element set $\{0,1\}$ are called \emph{Boolean functions.} If $A = \{0,1\}$ and $B$ is an arbitrary set, not necessarily equal to $\{0,1\}$, we refer to such functions as \emph{pseudo-Boolean functions.}

\subsection{Essential arity and arity gap}
\label{ssec:eagap}

Let $f$ be an $n$-ary $B$-valued function on $A$. For $1 \leq i \leq n$, the $i$-th variable is said to be \emph{essential} in $f$, or $f$ is said to \emph{depend} on $x_i$, if there are elements $a_1, \ldots, a_n, b \in A$ such that
\[
f(a_1, \ldots, a_{i-1}, a_i, a_{i+1} \ldots, a_n) \neq f(a_1, \ldots, a_{i-1}, b, a_{i+1}, \ldots, a_n).
\]
The number of essential variables in $f$ is called the \emph{essential arity} of $f$, and it is denoted by $\ess f$. 

We say that a function $f \colon A^n \to B$ is obtained from $g \colon A^m \to B$ by \emph{simple variable substitution,} or $f$ is a \emph{simple minor} of $g$, if there is a mapping $\sigma \colon \{1, \ldots, m\} \to \{1, \ldots, n\}$ such that
\[
f = g(x_{\sigma(1)}^{(n)}, \ldots, x_{\sigma(m)}^{(n)}).
\]
If $\sigma$ is not injective, then we speak of \emph{identification of variables.} If $\sigma$ is not surjective, then we speak of \emph{addition of inessential variables.} If $\sigma$ is a bijection, then we speak of \emph{permutation of variables.}
For indices $i, j \in \{1, \ldots, n\}$, $i \neq j$, the function $\minor{f}{i}{j} \colon A^n \to B$ obtained from $f \colon A^n \to B$ by the simple variable substitution
\[
\minor{f}{i}{j} := f(x_1^{(n)}, \ldots, x_{i-1}^{(n)}, x_j^{(n)}, x_{i+1}^{(n)}, \ldots, x_n^{(n)})
\]
is called a \emph{variable identification minor} of $f$, obtained by identifying $x_i$ with $x_j$.

The simple minor relation constitutes a quasi-order $\leq$ on the set of all $B$-valued functions of several variables on $A$ which is given by the following rule: $f \leq g$ if and only if $f$ is obtained from $g$ by simple variable substitution. If $f \leq g$ and $g \leq f$, we say that $f$ and $g$ are \emph{equivalent,} denoted $f \equiv g$. If $f \leq g$ but $g \not\leq f$, we denote $f < g$. It can be easily observed that if $f \leq g$ then $\ess f \leq \ess g$, with equality if and only if $f \equiv g$. For background, extensions and variants of the simple minor relation, see, e.g., \cite{CP2,Couceiro,CP,EFHH,FH,Lehtonen,LS,Pippenger, Wang, Zverovich}.

For $f \colon A^n \to B$ with $\ess f \geq 2$, we denote
\[
\essl f := \max_{g < f} \ess g,
\]
and we define the \emph{arity gap} of $f$ by $\gap f := \ess f - \essl f$. It is easily observed that
\[
\gap f = \min_{i \neq j} (\ess f - \ess \minor{f}{i}{j}),
\]
where $i$ and $j$ range over the set of indices of essential variables of $f$.
Since the arity gap is defined in terms of essential variables and since every $B$-valued function on $A$ is equivalent to a function all of whose variables are essential, we will assume without loss of generality that the functions $f \colon A^n \to B$ whose arity gap we consider are essentially $n$-ary.

The arity gap of $f$ is clearly at least $1$, and it can be as large as the number of essential variables of $f$, as illustrated by the following example due to Salomaa \cite{Salomaa}. Let $A$ be a finite set with $k$ elements ($k \geq 2$), and consider the operation $f \colon A^k \to A$ defined as follows. Let $\mathbf{b} = (b_1, \dots, b_k) \in A^k$ be a fixed $k$-tuple such that $b_i \neq b_j$ whenever $i \neq j$, and let $c$ and $d$ be distinct elements of $A$. We let
\[
f(\mathbf{a}) =
\begin{cases}
c, & \text{if $\mathbf{a} = \mathbf{b}$,} \\
d, & \text{otherwise.}
\end{cases}
\]
It is clear that $f$ depends on all of its $k$ variables, but whenever any pair of its variables is identified, the resulting function is a constant function, having no essential variables. Therefore, the arity gap of $f$ is $k$.
Thus, in order to classify all $B$-valued functions on $A$ according to their arity gap, we need to impose certain conditions on the functions being considered.  
One of such approaches requires that the functions are finite, i.e., that $A$ is finite.  
 
Partial results concerning the arity gap of finite functions were provided in  \cite{BK, GK, Willard}, and a general classification of finite functions according to their arity was given in \cite{CL2}. In the next subsection, we present such a characterization in the particular case of functions with Boolean variables.

\subsection{The arity gap of functions with Boolean variables}
\label{ssec:gapbf}

It is well-known that every Boolean function $f \colon \{0, 1\}^n \to \{0, 1\}$ is represented by a unique multilinear polynomial over the two-element field, the so-called Zhegalkin (or Reed--Muller) polynomial of $f$ \cite{Muller, Reed, Zhegalkin}.

\begin{theorem}[{\cite{CL}}]\label{BooleanGap}
Let $f \colon \{0, 1\}^n \to \{0, 1\}$ be a Boolean function with at least two essential variables. Then the arity gap of $f$ is $2$ if and only if $f$ is equivalent to one of the following functions:
\begin{enumerate}
\item $x_1 + x_2 + \dots + x_m + c$ for some $m \geq 2$,
\item $x_1 x_2 + x_1 + c$,
\item $x_1 x_2 + x_1 x_3 + x_2 x_3 + c$,
\item $x_1 x_2 + x_1 x_3 + x_2 x_3 + x_1 + x_2 + c$,
\end{enumerate}
where $c \in \{0,1\}$. Otherwise the arity gap of $f$ is $1$.
\end{theorem}

This complete classification of Boolean functions according to their arity gap in turn leads to the classification of pseudo-Boolean functions.

\begin{theorem}[{\cite{CL2}}]
\label{gappseudoB}
For a pseudo-Boolean function $f \colon \{0,1\}^n \to B$, $n \geq 2$, which depends on all of its variables, $\gap f = 2$ if and only if $f$ satisfies one of the following conditions:
\begin{enumerate}
\item $n = 2$ and $f$ is a nonconstant function satisfying $f(0,0) = f(1,1)$,
\item $f = g \circ h$, where $g \colon \{0, 1\} \to B$ is injective and $h \colon \{0, 1\}^n \to \{0, 1\}$ is a Boolean function with $\gap h = 2$, as listed above.
\end{enumerate}
Otherwise $\gap f = 1$.
\end{theorem}

\section{The arity gap of lattice polynomial functions}
\label{sec:gap}

In this section we study the arity gap of certain lattice functions, namely, the so-called polynomial functions.
Throughout this section,  let $L$ denote an arbitrary bounded distributive lattice (possibly infinite) with
lattice operations $\wedge$ and $\vee$, and with least and greatest elements $0$ and $1$, respectively. 
We denote by $\leqslant$ the associated lattice order.
For general background in lattice theory we refer the
reader to, e.g., Davey and Priestley \cite{DavPri02},
 Gr\"atzer~\cite{Grae03} and Rudeanu~\cite{Rud01}.

\subsection{Polynomial functions over bounded distributive lattices}
\label{ssec:polf}

By a (\emph{lattice}) \emph{polynomial function} we mean a map  $f \colon L^n \to L$ which can be obtained by
composition of the binary operations $\wedge$ and $\vee$, the projections, and the constant functions. If
constant functions are not used, then such polynomial functions are usually referred to as \emph{term functions}.
Note that polynomial functions are order-preserving, that is, 
$$
f(\vect{a})\leqslant f(\vect{b}) \quad \mbox{whenever} \quad  \vect{a}\leqslant \vect{b}.
$$
Moreover, the class of polynomial functions on a lattice $L$ is closed under formation of simple minors.

As shown by Goodstein  \cite{Goo67}, polynomial functions over bounded distributive lattices have very neat normal form representations.
Let $[n] := \{1, \ldots , n\}$. A lattice expression of the form
\begin{equation}
\bigvee_{I \subseteq [n]} a_I \bigwedge_{i \in I} x_i,
\label{eq:DNF}
\end{equation}
where the \emph{coefficients} $a_I$ ($I \subseteq [n]$) are elements of $L$, is said to be in \emph{disjunctive normal form} (\emph{DNF}).
A function $f\colon L^n \to L$ has a DNF representation, if there is an expression of the form \eqref{eq:DNF} which explicitly specifies $f$.
For instance, the \emph{median function} is represented in DNF by
\[
\median(x_1,x_2,x_3) = (x_1 \wedge x_2) \vee (x_2 \wedge x_3) \vee (x_3 \wedge x_1).
\]

Let $2^{[n]}$ denote the set of all subsets of $[n]$. For $I \subseteq [n]$, let $\vect{e}_I$ be the
\emph{characteristic vector} of $I$, i.e., the $n$-tuple in $L^n$ whose $i$-th component is $1$ if $i \in I$, and 0 otherwise.
Note that the mapping $\alpha \colon 2^{[n]}\to \{0,1\}^n$ given by $\alpha(I)=\vect{e}_I$, for every $I\in 2^{[n]}$, is an order-isomorphism. 

\begin{proposition}[{Goodstein \cite{Goo67}}]\label{prop:DNF(f)}
Let $L$ be a bounded distributive lattice. A function $f \colon L^n \to L$ is a polynomial function if and only if
\begin{equation}
f(x_1, \ldots ,x_n) = \bigvee_{I \subseteq [n]} \big( f(\vect{e}_I) \wedge \bigwedge_{i \in I} x_i \big).
\label{eq:Good}
\end{equation}
\end{proposition}

\begin{remark}
\label{rem:Good}
Observe that, by Proposition \ref{prop:DNF(f)}, every polynomial function $f \colon L^n \to L$ is uniquely determined by its restriction to $ \{0, 1\}^n$.
\end{remark}

\begin{remark}
Since every lattice polynomial function is order-preserving, we have that the coefficients in \eqref{eq:Good} are monotone increasing, i.e., $f(\vect{e}_I) \leqslant f(\vect{e}_J)$ whenever $I \subseteq J$. Moreover, a function $f \colon \{0,1\}^n \to L$ can be extended to a polynomial function over $L$ if and only if it is order-preserving.
\end{remark}

\subsection{Classification of lattice polynomial functions in terms of arity gap}
\label{ssec:gapl}

In this subsection, we will make use of Theorem \ref{gappseudoB} to obtain a complete classification of lattice polynomial functions in terms of arity gap.
We shall need the following auxiliary results.

\begin{proposition}
\label{prop:essential}
Let $L$ be a bounded distributive lattice and let $f \colon L^n \to L$ be a polynomial function.  Then for each $j \in [n]$, $x_j$ is essential in $f$ if and only if there exists a set $J \subseteq [n] \setminus \{j\}$ such that $ f(\vect{e}_J)< f(\vect{e}_{J \cup \{j\}})$.
\end{proposition}

\begin{proof}
Clearly, the condition is sufficient. To show that it is also necessary, assume that for all $J \subseteq [n] \setminus \{j\}$, we have that $ f(\vect{e}_J)= f(\vect{e}_{J \cup \{j\}})$. 
Consider the DNF representation of $f$ as given by equation \eqref{eq:Good}.
Then for every $J \subseteq [n] \setminus \{j\}$, the term $f(\vect{e}_J \cup \{j\})\wedge \bigwedge_{i\in J \cup \{j\}} x_i$ is absorbed by the term 
$f(\vect{e}_J )\wedge \bigwedge_{i\in J } x_i$. Hence, $f$ has a representation with no occurrence of $x_j$, and it is thus clear that $f$ does not depend on $x_j$.
\end{proof}

Using the description of the essential variables of polynomial functions given in Proposition \ref{prop:essential}, every polynomial function $f \colon L^n \to L$ has the same essential arity as the restriction of $f$  to $\{0, 1\}^n$.

\begin{proposition}\label{prop:ess2}
Let $L$ be a bounded distributive lattice, $f \colon L^n \to L$ a polynomial function and set $f' := f|_{\{0, 1\}^n}$. Then $x_j$ is essential in $f$ if and only if $x_j$ is essential in $f'$.
\end{proposition}

\begin{proof}
If $x_j$ is essential in $f'$, then there exist elements $a_1, \dots, a_n, a'_j \in \{0, 1\}$ such that
$f'(a_1, \dots, a_n) \neq f'(a_1, \dots, a_{j-1}, a'_j, a_{j+1}, \dots, a_n).$ 
Since $f(\vect{a}) = f'(\vect{a})$ for all $\vect{a} \in \{0, 1\}^n$, we have that $x_j$ is essential in $f$.

If $x_j$ is essential in $f$, then, by Proposition~\ref{prop:essential}, there exists a set $J \subseteq [n] \setminus \{j\}$ such that 
$ f(\vect{e}_J)< f(\vect{e}_{J \cup \{j\}})$. Since $f' = f|_{\{0, 1\}^n}$, we conclude that $x_j$ is essential in $f'$.
\end{proof}

Note that from Proposition \ref{prop:ess2}, it follows that the arity gap of a polynomial function $f$ coincides with the arity gap of its restriction $f|_{\{0, 1\}^n}$. In other words, the classification of polynomial functions with respect to their arity gap, reduces to that of pseudo-Boolean functions.
From Theorem \ref{gappseudoB}, we obtain the following explicit classification of polynomial functions.

\begin{theorem}\label{GapPol}
Let $L$ be a bounded distributive lattice and let $f \colon L^n \to L$ be a polynomial function with at least two essential variables. Then $\gap f = 2$ if and only if there are elements $a, b \in L$ such that $a < b$ and
\begin{equation}
\label{eq:gap2}
f \equiv (a \vee \median) \wedge b = a \vee ({\median} \wedge b).
\end{equation}
Otherwise, $\gap f = 1$.
\end{theorem}

\begin{proof}
As observed above, the arity gap of a polynomial function $f \colon L^n \to L$ is completely determined by its restriction $f' := f|_{\{0,1\}^n}$. By Theorem~\ref{gappseudoB}, the arity gap of $f'$ is either $1$ or $2$. If $f$ is of the form \eqref{eq:gap2}, then we clearly have that $\gap f = 2$. For the converse, assume that $\gap f = 2$. Since $f'$ is order-preserving, then using Theorem~\ref{gappseudoB} and ruling out the functions that are not order-preserving, we have that
\[
f' = g \circ h,
\]
where $g \colon \{0, 1\} \to L$ is injective and order-preserving and $h \equiv \median|_{\{0,1\}^3}$. Since the median function is idempotent and
\begin{align*}
f(\mathbf{0}) &= f'(\mathbf{0}) =: a, \\
f(\mathbf{1}) &= f'(\mathbf{1}) =: b,
\end{align*}
we have $g(0) = a$ and $g(1) = b$, and since $f$ is not a constant function, we have $a < b$. By Remark~\ref{rem:Good},
\[
f \equiv (a \vee \median) \wedge b = a \vee ({\median} \wedge b).
\]
\end{proof}


\begin{thebibliography}{99}
\bibitem{BK}
\textsc{J. Berman,} \textsc{A. Kisielewicz,}
On the number of operations in a clone,
\textit{Proc.\ Amer.\ Math.\ Soc.} \textbf{122} (1994) 359--369.

\bibitem{CP2}
\textsc{M. Bouaziz,} \textsc{M. Couceiro,} \textsc{M. Pouzet,}
Join-irreducible Boolean functions,
arXiv:\linebreak[0]0903.3848.

\bibitem{Cimev79}
\textsc{K.~N. \v{C}imev,}
On some properties of functions,
in: B. Cs\'ak\'any, I. Rosenberg (eds.)\ 
\textit{Finite Algebra and Multiple-Valued Logic,}
Abstracts of lectures of the colloquium on finite algebra and multiple-valued logic (Szeged, 1979),
North-Holland, 1981, pp.\ 38--40.

\bibitem{Cimev86}
\textsc{K.~N. \v{C}imev,}
\textit{Separable Sets of Arguments of Functions,}
Studies 180/1986,
Computer and Automation Institute, Hungarian Academy of Sciences,
Budapest, 1986.

\bibitem{Couceiro}
\textsc{M. Couceiro,}
On the lattice of equational classes of Boolean functions and its closed intervals,
\textit{J. Mult.-Valued Logic Soft Comput.} \textbf{18} (2008) 81--104.

\bibitem{CL}
\textsc{M. Couceiro,} \textsc{E. Lehtonen,}
On the effect of variable identification on the essential arity of functions on finite sets,
\textit{Int.\ J. Found.\ Comput.\ Sci.} \textbf{18} (2007) 975--986.

\bibitem{CL2}
\textsc{M. Couceiro,} \textsc{E. Lehtonen,}
Generalizations of \'{S}wierczkowski's lemma and the arity gap of finite functions,
\textit{Discrete Math.} \textbf{309} (2009)  5905--5912.

\bibitem{CP}
\textsc{M. Couceiro,} \textsc{M. Pouzet,}
On a quasi-ordering on Boolean functions,
\textit{Theoret.\ Comput.\ Sci.}\ \textbf{396} (2008) 71--87.

\bibitem{Davies}
\textsc{R.~O. Davies,}
Two theorems on essential variables,
\textit{J. London Math.\ Soc.}\ \textbf{41} (1966) 333--335.

\bibitem{DK}
\textsc{K. Denecke,} \textsc{J. Koppitz,}
Essential variables in hypersubstitutions,
\textit{Algebra Universalis} \textbf{46} (2001) 443--454.

\bibitem{DavPri02}
\textsc{B.~A. Davey,} \textsc{H.~A. Priestley,}
\textit{Introduction to Lattices and Order,}
2nd ed.,
Cambridge University Press, New York, 2002.

\bibitem{EFHH}
\textsc{O. Ekin,} \textsc{S. Foldes,} \textsc{P.~L. Hammer,} \textsc{L. Hellerstein,}
Equational characterizations of Boolean function classes,
\textit{Discrete Math.} \textbf{211} (2000) 27--51.

\bibitem{FH}
\textsc{A. Feigelson,} \textsc{L. Hellerstein,}
The forbidden projections of unate functions,
\textit{Discrete Appl.\ Math.} \textbf{77} (1997) 221--236.

\bibitem{Goo67}
\textsc{R.~L. Goodstein,}
The solution of equations in a lattice,
\textit{Proc.\ Roy.\ Soc.\ Edinburgh Sect.\ A} \textbf{67} (1965/1967) 231--242.

\bibitem{Grae03}
\textsc{G. Gr\"atzer,}
\textit{General Lattice Theory,}
2nd ed.,
Birkh\"auser Verlag, Berlin, 2003.

\bibitem{GK}
\textsc{G. Gr\"atzer,} \textsc{A. Kisielewicz,}
A survey of some open problems on $p_n$-sequences and free spectra of algebras and varieties,
in: A. Romanowska, J.~D.~H. Smith (eds.)\
\textit{Universal Algebra and Quasigroup Theory,}
Heldermann, Berlin, 1992, pp.\ 57--88.

\bibitem{Lehtonen}
\textsc{E. Lehtonen,}
Descending chains and antichains of the unary, linear, and monotone subfunction relations,
\textit{Order} \textbf{23} (2006) 129--142.

\bibitem{LS}
\textsc{E. Lehtonen,} \textsc{\'A. Szendrei,}
Equivalence of operations with respect to discriminator clones,
\textit{Discrete Math.}\ \textbf{309} (2009) 673--685.

\bibitem{Lupanov}
\textsc{O.~B. Lupanov,}
On a class of schemes of functional elements,
\textit{Problemy Kibernetiki} (1962) 61--114 (in Russian).

\bibitem{Muller}
\textsc{D.~E. Muller,}
Application of Boolean algebra to switching circuit design and to error correction,
\textit{IRE Trans.\ Electron.\ Comput.}\ \textbf{3}(3) (1954) 6--12.

\bibitem{Pippenger}
\textsc{N. Pippenger,}
Galois theory for minors of finite functions,
\textit{Discrete Math.}\ \textbf{254} (2002) 405--419.

\bibitem{Reed}
\textsc{I.~S. Reed,}
A class of multiple-error-correcting codes and the decoding scheme,
\textit{IRE Trans.\ Inf.\ Theory} \textbf{4}(4) (1954) 38--49.

\bibitem{Rud01}
\textsc{S. Rudeanu,}
\textit{Lattice Functions and Equations,}
Discrete Mathematics and Theoretical Computer Science Series,
Springer-Verlag, London, 2001.

\bibitem{Salomaa}
\textsc{A. Salomaa,}
On essential variables of functions, especially in the algebra of logic,
\textit{Ann.\ Acad.\ Sci.\ Fenn.\ Ser.\ A I.\ Math.}\ \textbf{339} (1963) 3--11.

\bibitem{Solovjev}
\textsc{N.~A. Solovjev,}
On the question of the essential dependence of functions of the algebra of logic,
\textit{Problemy Kibernetiki} \textbf{9} (1963) 333--335 (in Russian).

\bibitem{Wang}
\textsc{C. Wang,}
Boolean minors,
\textit{Discrete Math.}\ \textbf{141} (1991) 237--258.

\bibitem{Willard}
\textsc{R. Willard,}
Essential arities of term operations in finite algebras,
\textit{Discrete Math.}\ \textbf{149} (1996) 239--259.

\bibitem{Yablonski}
\textsc{S.~V. Yablonski,}
Functional constructions in a $k$-valued logic,
\textit{Tr.\ Mat.\ Inst.\ Steklova} \textbf{51} (1958) 5--142 (in Russian).

\bibitem{Zhegalkin}
\textsc{I.~I. Zhegalkin,}
On the calculation of propositions in symbolic logic,
\textit{Mat.\ Sb.}\ \textbf{34} (1927) 9--28 (in Russian).

\bibitem{Zverovich}
\textsc{I.~E. Zverovich,}
Characterizations of closed classes of Boolean functions in terms of forbidden subfunctions and Post classes,
\textit{Discrete Appl.\ Math.}\ \textbf{149} (2005) 200--218.
\end{thebibliography}
\end{document}